\documentclass[12pt]{amsart}
\usepackage{amsmath}
\usepackage{amsfonts}
\usepackage{amssymb}
\usepackage{amsthm}
\usepackage{graphicx}
\usepackage[all]{xy}

\makeatletter
\@namedef{subjclassname@2010}{%
\textup{2010} Mathematics Subject Classification}
\makeatother

\theoremstyle{plain}
\newtheorem{thm}{Theorem}[section]
\newtheorem{cor}[thm]{Corollary}
\newtheorem{lem}[thm]{Lemma}
\newtheorem{prop}[thm]{Proposition}

\theoremstyle{definition}
\newtheorem{de}[thm]{Definition}
\newtheorem{rem}[thm]{Remark}

\newcommand{\Spec}{\operatorname{Spec}}

\newcommand{\Res}{\operatorname{Res}}
\newcommand{\Hom}{\operatorname{Hom}}
\newcommand{\Rat}{\operatorname{Rat}}

\newcommand{\SL}{\operatorname{SL}}
\newcommand{\PGL}{\operatorname{PGL}}

\newcommand{\Aut}{\operatorname{Aut}}

\newcommand{\M}{\operatorname{M}}

\newcommand{\mc}[1]{\mathcal{#1}}
\newcommand{\mb}[1]{\mathbb{#1}}

\newcommand{\ti}[1]{\textit{#1}}
\newcommand{\tr}[1]{\textrm{#1}}
\newcommand{\tb}[1]{\textbf{#1}}
\newcommand{\ul}[1]{\underline{#1}}

\newcommand{\bd}{\begin{de}}
\newcommand{\ed}{\end{de}}
\newcommand{\bl}{\begin{lem}}
\newcommand{\el}{\end{lem}}
\newcommand{\bp}{\begin{prop}}
\newcommand{\ep}{\end{prop}}
\newcommand{\bt}{\begin{thm}}
\newcommand{\et}{\end{thm}}
\newcommand{\bc}{\begin{cor}}
\newcommand{\ec}{\end{cor}}
\newcommand{\bpf}{\begin{proof}}
\newcommand{\epf}{\end{proof}}

\newcommand{\beq}{\begin{equation}}
\newcommand{\eeq}{\end{equation}}
\newcommand{\beqs}{\begin{equation*}}
\newcommand{\eeqs}{\end{equation*}}
\newcommand{\ben}{\begin{enumerate}}
\newcommand{\een}{\end{enumerate}}
\newcommand{\bit}{\begin{itemize}}
\newcommand{\eit}{\end{itemize}}

\begin{document}

\baselineskip=17 pt

\title{The Moduli Space of Totally Marked \\ Degree Two Rational Maps}

\author{Anupam Bhatnagar}
\address{Department of Mathematics, Borough of Manhattan
Community College, The City University of New York; 199 Chambers Street, New York, NY 10007 U.S.A.}
\email{anupambhatnagar@gmail.com}

\date{\today}


\begin{abstract}
A rational map $\phi: \mb{P}^1 \to \mb{P}^1$ along with an ordered list of fixed and critical points
is called a totally marked rational map. The space of totally marked degree two rational maps, 
$\Rat^{tm}_2$ can be parametrized by an affine open subset of $(\mb{P}^1)^5$. 
We consider the natural action of $\SL_2$ on $\Rat^{tm}_2$ induced from the action of $\SL_2$ on
$(\mb{P}^1)^5$ and prove that the quotient space $ \Rat^{tm}_2/\SL_2$ exists as a scheme. The quotient
is isomorphic to a Del Pezzo surface with the isomorphism being defined over $\mb{Z}[1/2]$.
\end{abstract}



\maketitle


\section{Introduction}
A rational map $\phi: \mb{P}^1 \to \mb{P}^1$  of degree two over a field $k$ is given by a pair 
of homogeneous polynomials
\begin{align*}
\phi = [\phi_0, \phi_1] = [aX^2 + bXY + c Y^2, dX^2 + eXY+fY^2]
\end{align*}
such that $\phi_0, \phi_1$ have no common roots. In non-homogeneous form $\phi$ may be expressed as
\begin{align*}
\phi(z) = \frac{az^2+bz+c}{dz^2+ez+f}
\end{align*}
Let $\phi_0(z) = az^2+bz+c$ and $\phi_1(z)=dz^2+ez+f$. We define $\Res(\phi)$, the resultant of $\phi$ 
as the product 
$$ \prod_{(\alpha,\beta): \phi_0(\alpha)=\phi_1(\beta)=0} (\alpha - \beta)$$
 The condition that $\phi_0,\phi_1$ have no
common roots is equivalent to $\Res(\phi) \neq 0$.

Let $\Rat_2$ denote the space of degree two rational maps $\mb{P}^1 \to \mb{P}^1$.
The special linear group $\SL_2$ acts via conjugation on $\Rat_2$ i.e. for $f\in \SL_2, \phi \in \Rat_2,
f \cdot \phi = f \circ \phi \circ f^{-1}$. The moduli space $\Rat_2/\SL_2$, denoted $\M_2$ 
arises naturally in the study of dynamical systems on $\mb{P}^1$. Over the complex numbers
Milnor proved in \cite{Mi1} that $\Rat_2(\mb{C})/\SL_2(\mb{C})$ is biholomorphic to $\mb{C}^2$. 
This fact was generalized by Silverman in \cite{S}; he showed $\M_2$ is an affine
integral scheme over $\mb{Z}$ and is isomorphic to $\mb{A}^2_{\mb{Z}}$. 

Inspired by Milnor \cite{Mi2} we consider a rational map along with an ordered list of its fixed and 
critical points. Since a rational map of degree two is completely determined by its fixed and critical
points, we dispose of the map and focus on ordered lists of fixed and critical points. We refer to this 
as the space of totally marked degree two rational maps, $\Rat^{tm}_2$. It can 
be viewed as an affine open subvariety of $(\mb{P}^1)^5$.

Let $p_1,p_2,p_3,q_1,q_2$ be an ordered list of fixed points and critical points of some degree two 
rational map. The natural action of the special linear group $\SL_2$ on $(\mb{P}^1)^5$ induces an action
on $\Rat^{tm}_2$. In this article we analyze the quotient $\Rat^{tm}_2/\SL_2$ and prove:

\bt \label{main_thm_intro}
Let $\Rat^{tm}_2$ denote the space of totally marked degree two rational maps. Consider the following 
action of $\SL_2$ on $\Rat^{tm}_2$
$$f \cdot(p_1,p_2,p_3, q_1,q_2) =(f(p_1), f(p_2), f(p_3), f(q_1), f(q_2))$$
Then the moduli space $\Rat^{tm}_2/\SL_2$ is isomorphic to a Del Pezzo surface and the isomorphism is 
defined over $\mb{Z}[1/2]$.
\et

Recall that a cubic in $\mb{P}^3$ is a Del Pezzo surface. We give the explicit equation of the surface in \S5.
The above theorem generalizes a similar result by Milnor \cite{Mi2} over $\mb{C}$. The two most significant 
facts which allow us to prove the theorem above are:
\ben
\item[(a)] The fixed points and critical points of a degree two rational map determine the map completely.
\item [(b)] The three cross ratio's formed by selecting both critical points and selecting two out of the three fixed points at a time (see Definition \ref{rat2tm}) are $\SL_2$ invariant functions on $\Rat^{tm}_2$.
\een

Observe that for $\phi: \mb{P}^1 \to \mb{P}^1, z \mapsto z^2$ each point of $\mb{P}^1$ is a critical point in 
characteristic two. Thus the notion of a totally marked rational map is not well defined in characteristic two, 
so the isomorphism in theorem above cannot be defined over $\mb{Z}$.  

The moduli space of totally marked degree two rational maps, $\M^{tm}_2$ is a 12-to-1 cover of $\M_2$. 
Indeed, the map $\M^{tm}_2 \to \M_2$ factors through the moduli space of fixed point marked degree 
two rational maps, $\M^{fm}_2$.  $\M^{fm}_2$ is a 6-to-1 cover of $\M_2$ and $\M^{tm}_2$ is a double cover of $ \M^{fm}_2$.

It is natural to ask about the structure of the quotient $\M^{tm}_d: =\Rat^{tm}_d/\SL_2$. To answer 
this, we need analogs of (a) and (b) for $d>2$. As in the degree two case, $\M^{tm}_d$ will be a finite 
cover of $\M_d$ and studying $\M^{tm}_d$ is useful for finding equations defining $\M_d$.

In \S2 we prove some basic facts about degree two rational maps. In \S3 we decribe the moduli scheme of 
totally marked degree two rational maps $\M^{tm}_2$, followed by the moduli functor for totally marked 
degree two rational maps, $\ul{\M}^{tm}_2$ in \S4. We prove that the moduli scheme $\M^{tm}_2$ 
is a coarse moduli scheme for the functor $\ul{\M}^{tm}_2$. Finally in \S 5 we prove our main 
result.

\tb{Notation/Conventions}: Thoughout this article we fix $k$ to be a field of characteristic 
different from two. We denote the fixed points by $p_1, p_2, p_3$ and critical points by
$q_1, q_2$.

\section{Preliminaries}

\bl \label{distinct critical points}
Let $\phi: \mb{P}^1 \to \mb{P}^1$ be a rational map of degree two defined over 
$\mb{Z}[1/2]$ such that the resultant of $\phi$, $\Res(\phi)$ is nonzero. Then $\phi$ has 
two distinct critical points.
\el
\bpf
Let $$\phi(z) = \frac{az^2+bz+c}{dz^2+ez+f}$$
and denote the fixed points of $\phi$ by $p_1,p_2,p_3$. We split the proof in two cases.

Case 1: Suppose there is a fixed point of multiplicity three. Without loss of generality we may assume 
$p_1$ has multplicity three. Then
$$ \phi (z) =  \frac{az^2+bz-p^3_1}{z^2+(a-3p_1)z+(b+3p_1^2)}$$
with an appropriate change of coordinates if $p_1 = \infty$.

Case 2: Suppse there is no fixed point with multiplicity three. Without loss of generality we may assume 
that $p_2 \neq p_3$. Applying a change of coordinates we let  $p_2=0, p_3=\infty$. Then 
$$ \phi (z) = \frac{az^2+bz}{ez+f}$$
In both cases it can be easily verified that the critical points are distinct. 
\epf

\bl \label{unique ratl map}
Let $\phi: \mb{P}^1 \to \mb{P}^1$ be a rational map of degree two defined over $\mb{Z}[1/2]$. Then 
$\phi$ is uniquely determined by its fixed points and critical points.
\el
\bpf
Let  $$\phi(z) = \frac{az^2+bz+c}{dz^2+ez+f}$$ and denote its fixed and critical points by 
$p_1, p_2, p_3$ and $q_1,q_2$ respectively. By the previous lemma we know $q_1 \neq q_2$,
so we may assume $q_1= 0$ and $q_2 = \infty$. Observe that 
\[
\begin{array}{rcl}
q_1 =0 \tr{ and } q_2 = \infty& \iff & ae-bd =0 \tr{ and } bf-ce=0\\
&\iff & \phi(z) = \phi(-z) \tr{ for all } z \in \mb{P}^1 \\
& \iff & b= e= 0
\end{array}
\]
Therefore, $$ \phi (z)= \frac{a z^2 + c}{d z^2 + f} $$
There is a fixed point at infinity if and only if $d=0$. The fixed points of $\phi$ are the roots of 
the equation $dz^3 -a z^2 + f z - c=0$, and they uniquely determine the point $(d:a:f:c)$ in 
$\mb{P}^3$. Thus the coefficients $a,c,d,f$ and hence the rational map $\phi$ is uniquely 
determined by its fixed points and critical points.
\epf

\section{The Moduli Scheme $\M^{tm}_2$}

For any vector $v= (v_1, \ldots, v_n) \in \mb{Z}^n$ we define a line bundle on 
$(\mb{P}^1)^n$ $$ L_{v} = \bigotimes_{i=1}^n \pi^*_i ( \mc{O}_{\mb{P}^1}(1)^{\otimes v_i}) $$ 
where $\pi_i: (\mb{P}^1)^n \to \mb{P}^1$ is the projection on the $i^{th}$ factor.

\bd \label{rat2tm}
Let $(\omega_1, \omega_2, \omega_3, \xi_1, \xi_2)$ be non-homogeneous coordinates
on $(\mb{P}^1)^5$. Fix the linearization $m=(1,1,1,2,2)$ on $(\mb{P}^1)^5$ and denote it by 
$(\mb{P}^1)^5(L_m)$. Let 
$$C:= \{ (\omega_1, \omega_2, \omega_3, \xi_1, \xi_2) \in (\mb{P}^1)^5(L_m) \;|\;
\xi_1 = \xi_2 \}$$ and $$R_i:= \bigg\{ (\omega_1, \omega_2, \omega_3, \xi_1, \xi_2) \in 
(\mb{P}^1)^5(L_m) \; \bigg|\; r_i:=
\frac{(\omega_j - \xi_1)(\omega_k - \xi_2)}{(\omega_j -\xi_2)(\omega_k -\xi_1)} =  -1 \bigg\} $$
where $(i,j,k)$ is any permutation of (1,2,3). We define the \ti{space of totally marked degree two rational maps} as $$ \Rat^{tm}_2 := (\mb{P}^1)^5 (L_m)\setminus \{C \cup 
R_1 \cup R_2 \cup R_3\} $$
\ed

A generic element of $\Rat^{tm}_2$ is an ordered set of fixed points and critical 
points of a degree two rational map. Observe that if $\xi_1=0$ and $\xi_2= \infty$, 
then $\omega_i \neq -\omega_j$ for $i \neq j$. The automorphism group of $\mb{P}^1$, 
$\PGL_2$ acts on each coordinate of $\Rat^{tm}_2$. For technical reasons we consider
the action of $\SL_2$ instead of $\PGL_2$.

\bd
Two elements $\{p_1, p_2, p_3, q_1, q_2\}$ and $\{p'_1, p'_2, p'_3, q'_1,q'_2 \}$ 
of $\Rat^{tm}_2$ are said to be $\SL_2$ equivalent if there exists 
$f\in \SL_2$ such that $f(p_i) = p'_i, f(q_i) = q'_i$. The quotient $\Rat^{tm}_2 /\SL_2$ is called the 
\ti{moduli space of totally marked degree two rational maps} and denoted by $\M^{tm}_2$. 
\ed

A priori, for an algebraically closed field $k$ the quotient $\M^{tm}_2(k)$ exists 
as a set. We shall show that this set is isomorphic to a Del Pezzo surface, whenever
char($k$) $\neq 2$. We now describe the set of stable and semistable points of 
projective space. This is well known, we recall it here for the reader's convenience.
The set of stable and semistable points of a scheme (say $V$) are denoted by $V^s$ and $V^{ss}$ 
respectively.

\bt
Let $P = (x_1, \ldots, x_m) \in (\mb{P}^r)^m$ and $v= (v_1, \ldots, v_m) \in \mb{Z}^m$. Then 
$$ P \in ((\mb{P}^r)^m)^{ss}(L_v)\;\; (resp. \;\; P \in ((\mb{P}^r)^m)^{s}(L_v) )$$
if and only if for every proper linear subspace $W$ of $\mb{P}^r$
$$ \sum_{i,x_i \in W} v_i \leq \frac{(\dim W +1)}{n+1} \bigg(\sum_{i=1}^m v_i \bigg) $$
(resp. the strict inequality holds). 
\et

\bpf
\cite{D}, p. 172.
\epf

\bc \label{semistable}
$$  ((\mb{P}^r)^m)^{ss}(L_v) \neq \emptyset \iff \forall i = 1, \ldots, m,\;\; (r+1) v_i \leq 
\sum_{i=1}^m v_i $$ 
$$  ((\mb{P}^r)^m)^{s}(L_v) \neq \emptyset \iff \forall i = 1, \ldots, m,\;\; (r+1) v_i <
\sum_{i=1}^m v_i $$ 
\ec

\bpf
\cite{D}, p. 172.
\epf
Using Corollary \ref{semistable} with the  linearization $m=(1,1,1,2,2)$ we have,
\beq \label{*}
\Rat^{tm}_2 \subset ((\mb{P}^1)^5)^s(L_m) = ((\mb{P}^1)^5)^{ss}(L_m)
\eeq
The equality in (1) follows from Corollary \ref{semistable} and the inclusion $\Rat^{tm}_2 
\subset ((\mb{P}^1)^5)^s(L_m)$ follows by observing that $(1,-1,2,0,\infty)\in 
((\mb{P}^1)^5)^s(L_m)$ but $(1,-1,2,0, \infty) \notin \Rat^{tm}_2$ since 
$\omega_1 = -\omega_2$. The choice of linearization $(1,1,1,2,2)$ is not arbitrary. 
If we use the linearization $(1,1,1,1,1)$, then by Corollary \ref{semistable} it can be 
verified that $(1,1,1,0,\infty) \notin ((\mb{P}^1)^5)^{ss}(L_m)$. The rational map 
$(3z^2+1)/(z^2+3)$ has a triple fixed point at 1 and critical points at $0$ and $\infty$. 

\bt \label{quotient}
Using the notation above and the linearization $m=(1,1,1,2,2)$ we have:
\ben
\item[(a)]
The space of totally marked degree two rational maps, $\Rat^{tm}_2$ is an 
$\SL_2$-invariant open subset of the stable locus $((\mb{P}^1)^5)^s(L_m)$ in 
$(\mb{P}^1)^5(L_m)$. Hence, the geometric quotient 
$\M^{tm}_2 = \Rat^{tm}_2/\SL_2$ 
exists as a scheme over $\mb{Z}[1/2]$.

\item[(b)]
The geometric quotient $(\M^{tm}_2)^s = ((\mb{P}^1)^5)^s(L_m)/\SL_2$ and the 
categorical quotient  $(\M^{tm}_2)^{ss}=  ((\mb{P}^1)^5)^{ss}(L_m)/\SL_2$ exist
as schemes over $\mb{Z}[1/2]$ and are the same for the linearization $(1,1,1,2,2)$.

\item[(c)]
The schemes $\M^{tm}_2, (\M^{tm}_2)^s$ and $(\M^{tm}_2)^{ss}$ are connected,
integral, normal and of finite type over $\mb{Z}[1/2]$. Moreover, $\M^{tm}_2$
is affine over $\mb{Z}[1/2]$.
\een
\et

\bpf The proof follows from standard invariant theoretic results in \cite{Mu} and \cite{Se}.
\ben 
\item[(a)]
The inclusion $\Rat^{tm}_2 \subset ((\mb{P}^1)^5)(L_m)$ follows from (\ref{*}). The action of $\SL_2$
on $(\mb{P}^1)^5(L_m)$ fixes the sets $R_i$  and $C$ defined in Definition \ref{rat2tm}. Hence, 
$\Rat^{tm}_2$ is an $\SL_2$-stable and $\SL_2$ invariant scheme, so the geometric quotient 
$\M^{tm}_2 = \Rat^{tm}_2/\SL_2$ exists. 
Over a field this a consequence of Mumford's construction of quotients \cite{Mu}, Chapter 1, and over 
$\mb{Z}[1/2]$ it follows by essentially the same methods using Seshadri's theorem that a reductive 
group scheme is geometrically reductive (see \cite{Mu} and \cite{Se}).

\item[(b)]
The existence of quotients follows from Mumford \cite{Mu} and Seshadri \cite{Se} and 
the equality $(\M^{tm}_2)^s = (\M^{tm}_2)^{ss}$ follows from Corollary \ref{semistable}.

\item[(c)]
The schemes $\Rat^{tm}_2, ((\mb{P}^1)^5)^s$ and $((\mb{P}^1)^5)^{ss}$ are open subschemes of 
$(\mb{P}^1)^5$ so they are connected, integral and normal. By \cite{Mu}, Section 2, Remark 2, we 
conclude that the respective quotients $\M^{tm}_2, (\M^{tm}_2)^s$ and $(\M^{tm}_2)^{ss}$ are
connected, integral and normal. The fact that $\M^{tm}_2$ is affine over $\mb{Z}[1/2]$ follows from 
\cite{Mu}, Theorem 1.1.

\een
\epf


\section{The Moduli Functor $\ul{\M}^{tm}_2$}

\bd 
The functor $\ul{\Rat}^{tm}_2$ of totally marked degree two rational maps is the functor 
$$ \ul{\Rat}^{tm}_2: (Sch/\mb{Z}[1/2]) \to (Sets)$$ defined by 
\[
\ul{\Rat}^{tm}_2(S) =
\left\{
\begin{array}{l}
\tr{separable $S$-morphisms } \phi: \mb{P}^1_S \to \mb{P}^1_S \tr{ satisfying }\phi^*\mc{O}(1) \cong \mc{O}(2), \\
\tr{sections } r_i \tr{ of }\mb{P}^1_S \to S, i=1,2,3 \tr{ satisfying } \phi \circ r_i = r_i, \tr{ and } \\
\tr{sections } s_j \tr{ of }\mb{P}^1_S \to S, j= 1,2 \tr{ satisfying } div(s_1)+div(s_2) = R_{\phi}
\end{array}
\right.
\]
where $R_\phi$ is the ramification divisor of $\phi$.
\ed

Observe that the sections $r_i$ correspond to the fixed points and the sections $s_j$ correspond to 
critical points. A $S$-point of $\ul{\Rat}^{tm}_2$ consists of a degree two rational map and five 
sections satisfying the above conditions. 

\bd We say two $S$-points of $\ul{\Rat}^{tm}_2$, say $(\phi, r_1,r_2, r_3, s_1, s_2)$ and 
$(\phi', r'_1,r'_2, r'_3, s'_1, s'_2)$ are equivalent if there exists $f \in \Aut(\mb{P}^1_S)$
such that $\phi \circ f = f \circ \phi', f(r_i) = r'_i,$ and $f(s_j) = s'_j$. We define 
the moduli functor $\ul{\M}^{tm}_2$ to be the quotient of $\ul{\Rat}^{tm}_2$ under the above
equivalence relation.
$$ \ul{\M}^{tm}_2 : (Sch/\mb{Z}[1/2]) \to (Sets) \hspace{10mm} S \mapsto \ul{\Rat}^{tm}_2(S)/ \sim $$
\ed

We now prove that the functor $\ul{\Rat}^{tm}_2$ is representable.

\bt
The scheme $\Rat^{tm}_2$ defined in Definition \ref{rat2tm} represents the functor $\ul{\Rat}^{tm}_2$. In particular, $\Rat^{tm}_2$ is a fine moduli space for $\ul{\Rat}^{tm}_2$.
\et

\bpf 
Let $S$ be an arbitary scheme and $(p_1,p_2,p_3, q_1,q_2) \in \Rat^{tm}_2(S)$. By Lemma 
\ref{unique ratl map} there exists a unique rational map $\phi: \mb{P}^1_S \to \mb{P}^1_S$ 
with fixed points $p_1,p_2,p_3$ and critical points $q_1,q_2$. Let $r_i, s_j$ be the sections 
of $\mb{P}^1_S \to S$ corresponding to the fixed and critical points. This gives a well 
defined map
\begin{align} \label{scheme to functor}
\Rat^{tm}_2(S)  &\to \ul{\Rat}^{tm}_2(S) \\
(p_1,p_2,p_3,q_1,q_2) &\mapsto (\phi, r_1,r_2r_3, s_1,s_2) \notag
\end{align}
The inverse 
\begin{align} \label {functor to scheme}
 \ul{\Rat}^{tm}_2(S) \to & \Rat^{tm}_2(S) \\
(\phi, r_1,r_2, r_3, s_1, s_2) \mapsto & (p_1,p_2,p_3,q_1,q_2) \notag
\end{align}
maps the sections $r_i, s_j$ to the corresponding fixed and critical points and forgets $\phi$.
Thus the scheme $\Rat^{tm}_2$ represents the functor $\ul{\Rat}^{tm}_2$.
\epf

We now show $\M^{tm}_2$ is a coarse moduli scheme for the functor $\ul{\M}^{tm}_2$. 

\bt
There is a natural map of functors $$\ul{\M}^{tm}_2 \to \Hom(-, \M^{tm}_2)$$
with the property that $\ul{\M}^{tm}_2(k) \cong \M^{tm}_2(k)$ for every algebraically
closed field k of characteristic $\neq 2$. 
\et

\bpf 

Let $S$ be a arbitrary scheme over $\mb{Z}[1/2], [\eta] \in \ul{\M}^{tm}_2(S)$ and $(\phi, r_1, r_2, r_3,s_1, s_2) 
\in \ul{\Rat}^{tm}_2(S)$ be a representative of $[\eta]$. Let $(p_1,p_2,p_3,q_1,q_2) \in \Rat^{tm}_2$
be the image of  $(\phi, r_1, r_2, r_3,s_1, s_2)$ along the map defined in (\ref{functor to scheme}). Taking 
the quotient of $(p_1,p_2,p_3,q_1,q_2)$ by $\SL_2$ we get the image of $[\eta]$ in $\Hom(-,\M^{tm}_2)$. 
The image is independent of the choice of the lifting of $[\eta]$ since we quotient by $\SL_2$. 

For any algebraically closed field $k$ of characteristic different than two, 
$$ \ul{\M}^{tm}_2(k) \cong \Rat^{tm}_2(k)/\PGL_2(k), \hspace{10 mm} 
\M^{tm}_2(k) \cong \Rat^{tm}_2(k)/\SL_2(k)$$
The map $\SL_2 \to \PGL_2$ is surjective hence the quotients are the same.
\epf


\section{Main Theorem}

Recall that $\Rat^{tm}_2 :=  (\mb{P}^1)^5(L_m) \setminus \{C \cup R_1 \cup R_2 \cup R_3\} $, where 
$$C:= \{ (\omega_1, \omega_2, \omega_3, \xi_1, \xi_2) \in (\mb{P}^1)^5(L_m) \;|\;
\xi_1 = \xi_2 \}$$ and $$R_i:= \bigg\{ (\omega_1, \omega_2, \omega_3, \xi_1, \xi_2) \in 
(\mb{P}^1)^5(L_m) \; \bigg|\; r_i:=
\frac{(\omega_j - \xi_1)(\omega_k - \xi_2)}{(\omega_j -\xi_2)(\omega_k -\xi_1)} =  -1 \bigg\} $$
where $(i,j,k)$ is any permutation of (1,2,3). We shall show that the cross ratios $r_i$ are 
$\SL_2$--invariant functions on the quotient space $\M^{tm}_2$ and they uniquely determine 
the conjugacy class.

\bp
The cross ratios $r_i$ are rational functions on 
$\Rat^{tm}_2$. Moreover, they are invariant under the action of $\SL_2$ on $\Rat^{tm}_2$. 
Thus they descend to give rational functions on $\M^{tm}_2$. 
\ep

\bpf 
Two elements $(p_1,p_2,p_3,q_1,q_2), (p'_1,p'_2,p'_3,q'_1,q'_2) \in \Rat^{tm}_2$ are $\SL_2$ equivalent if 
there exists $f\in \SL_2$ such that $f(p_i) = p'_i$ and $f(q_i) = q'_i$, where $p_i$ denote the fixed points and 
$q_i$ denote the critical points. Note that each cross ratio is determined 
by selecting two of the three fixed points and both critical points. If 
$$ r_1 = \frac{(p_2 - q_1)(p_3 - q_2)}{(p_2 - q_2)(p_3 -q_1)}$$
then the cross ratio determined by $f(p_2), f(p_3), f(q_1), f(q_2)$ is 
\begin{align} \label{cross ratio}
r'_1 = \frac{(f(p_2) - f(q_1))(f(p_3) - f(q_2 ))}{(f(p_2) - f(q_2))(f(p_3) - f(q_1))} = 
\frac{(p'_2 - q'_1)(p'_3 - q'_2)}{(p'_2 - q'_2)(p'_3 -q'_1)}
\end{align}

Claim: $r_1 = r'_1$

Proof of Claim: By Lemma \ref{distinct critical points} we may assume $q_1 =0, q_2 = \infty$,
hence $r_1= p_2/p_3$. 
Write $p_i = [p_i:1], \infty = [1:0]$ and  $0=[0:1]$. For 
$$ f=
\left[
\begin{array}{cc}
a & b  \\
c & d  
\end{array}
\right] \in \SL_2
$$
$f(p_i) = [\frac{ap_i+b}{cp_i+d}:1], f(\infty) = [a/c:1]$ and $f(0)= [b/d:1]$. Writing these in 
non-homogeneous form and substituting in (\ref{cross ratio}), we get $r'_1 = p_2/p_3$. Similarly for $r_2, r_3$. 
Since the cross ratio's are invariant under the $\SL_2$ action, they descend to give rational functions
on $\M^{tm}_2$.
\epf

Let $V= \Spec( \mb{Z}[1/2][x_1,x_2, x_3]/(x_1+x_2+x_3+x_1 x_2 x_3))$. 

\bp
The cross ratios form a complete conjugacy invariant i.e. they determine the conjugacy class in 
$\M^{tm}_2$ uniquely. 
\ep

\bpf 
We begin by defining a map from the scheme $V$ to the fixed point marked moduli space, denoted 
$\M^{fm}_2$  and then extending it to $\M^{tm}_2$. The fixed point marked moduli space is determined
by the multipliers at the three fixed points which we denote by $\mu_1, \mu_2, \mu_3$. Define a 
map from $V$ to $\M^{fm}_2$ by setting 
$$ \mu_i = 1+x_j x_k$$
Observe that $\mu_1 + \mu_2 + \mu_3 = \mu_1\mu_2\mu_3 +2$. If $\omega_j \neq \omega_k$, then
we can put $\omega_j =0$ and $\omega_k= \infty$ and write the map as 
$$\phi(z) = \frac{z^2 + \mu_j z}{\mu_k z+1}$$
The critical points of $\phi$ are $$\xi_1 = \frac{-1 +  x_i}{\mu_k} \hspace{10 mm} 
\xi_2 = \frac{-1- x_i}{\mu_k}$$
and the cross ratio $r_i$ is given by $$ r_i = \frac{\xi_1}{\xi_2} = \frac{1-x_i}{1+x_i}$$
Conversely, given $r_i$ we can solve for $x_i = (1-r_i)/(1+r_i)$. This shows the cross ratios and hence the
conjugacy class in $\M^{tm}_2$ is completely determined by the coordinates $x_1,x_2,x_3$, yielding a
smooth map from $\M^{tm}_2$ to $V$.
\epf

Let $r_1,r_2, r_3$ be non-homogeneous coordinates on $(\mb{P}^1 \setminus \{-1\})^3$
and let $W$ be the subvariety cut out by the equation $r_1r_2r_3 -1$. We now prove that 
$\M^{tm}_2$ is isomorphic to $W$, where the isomorphism is defined over $\mb{Z}[1/2])$. 

\begin{rem}
The schemes $V$ and $W$ are isomorphic to each other. Using $x_1,x_2,x_3$ and 
$r_1,r_2,r_3$ as coordinates on $V$ and $W$ respectively define $\sigma: V \to W,
x_i \mapsto r_i = (1-x_i)/(1+x_i)$. It can be easily verified that $\sigma= \sigma^{-1}$. 
\end{rem}

\bt
The map \begin{align*}
\M^{tm}_2 & \to W 
 \end{align*}
is an isomorphism of schemes defined over $\mb{Z}[1/2]$.
\et

\bpf
Let $(\omega_1, \omega_2, \omega_3, \xi_1, \xi_2)$ be any element of 
$\M^{tm}_2$. The map from $\M^{tm}_2 \to W$ is given by 
$$r_i = \frac{(\omega_j - \xi_1)(\omega_k - \xi_2)}{(\omega_j -\xi_2)(\omega_k -\xi_1)}$$
where $(i,j,k)$ is any permutation of $(1,2,3)$.
We now construct the inverse. Without loss of generality we may assume that one of 
$\omega_1, \omega_2, \omega_3$ is finite and nonzero, and $\omega_1 =1, 
\omega_2 = 1/r_3$, $\omega_3 = r_2, \xi_1=0, \xi_2 = \infty$.
Since $\xi_1 =0, \xi_2 = \infty$ we have $\phi(z)= \frac{az^2+b}{cz^2+d}$. We shall determine the
coefficients of $\phi$ explicitly. The image for these values of $\omega_1, \omega_2, \omega_3, \xi_1, \xi_2$
in $(\mb{P}^1\setminus \{-1\})^3$ is 
the complement of the curves $(r_2 = \infty, r_3=0)$ and $(r_2=0, r_3 = \infty)$. We denote the image in 
$(\mb{P}^1 \setminus \{-1\})^3$ by  $U$. 
For $\omega_1, \omega_2, \omega_3, \xi_1, \xi_2$ as above, we have: 
$$a/c= -(1+r_2+ 1/r_3), b/c = r_2/r_3, d/c = r_2 + 1/r_3 + r_2/r_3$$
We now break $U$ into four subsets based on values of $r_2$ and $r_3$. 
\ben
\item[Case 1:] If $r_2 \neq \infty, r_3 \neq 0$, then let $c=1$ and we are done.
\item[Case 2:] If $r_2 \neq \infty, r_3 \neq \infty$, then let $c=r_3 \implies b = -r_2, a= -(1+r_3 + r_2r_3), d=1+r_2+r_2r_3$.
\item[Case 3:] If $r_2 \neq 0, r_3 \neq 0$, then let $b=-1/r_3 \implies a= -(1+ 1/r_2+ 1/r_2r_3), c = 1/r_2, d= 1+ 1/r_3 + 1/r_2r_3$.
\item[Case 4:] If $r_2\neq 0, r_3 \neq \infty$, then let $b =-1 \implies a = -(r_3 + 1/r_2 + r_3/r_2), c = r_3/r_2, d = 1+ r_3 + 1/r_2$.
\een
In each of the four cases $\phi(z) \in \M^{tm}_2$. On the intersections the maps agree not only 
in $\M^{tm}_2$ but also in $\Rat^{tm}_2$. So we can glue the four affine pieces together to obtain
a map $U \to \M^{tm}_2$. By symmetry we can assume that $\omega_2, \omega_3$ are finite and
nonzero as well. We get morphisms from three affine open pieces to $\M^{tm}_2$. The union of 
these three affine pieces is $W$, so we have three morphisms from $W$ to $\M^{tm}_2$. It remains to show
that the morphisms agree on the intersections.

On the first affine piece (i.e. $\omega_1 \neq 0, \infty$) we have $\omega_1 =1,
\omega_2 = 1/r_3, \omega_3 = r_2$, so $a/c= -(1+ r_2+1/r_3), b/c = -r_2/r_3,  
d/c = r_2 + 1/r_3 + r_2/r_3$.
On the second affine piece (i.e. $\omega_2 \neq 0, \infty$) we have $\omega_2 =1,
\omega_1 = 1/r_3, \omega_3 = r_2r_3$, so $a/c =-(1+r_3+r_2r_3), b/c = -r_2r^2_3, 
d/c = r_3(1+r_2+r_2r_3)$. Applying the transformation $z \mapsto r_3 \cdot z$ 
to the equation $\phi(z) =z$ we see the expressions for $a/c, b/c, d/c$ are the same.
On the third affine piece (i.e. $\omega_3 \neq 0, \infty$) if $r_2\neq \infty$, then let`
$c=1$ and if $r_2 \neq 0$, then let $c=1/r_2$. In either case $r_1,r_2,r_3$ determine
the same quadruple $(a,b,c,d)$ defining the same point in $\M^{tm}_2$, though not 
the same point in $\Rat^{tm}_2$.
\epf

A cubic in $\mb{P}^3$ is a Del Pezzo surface. In homogenoeous coordinates 
the surface $W$ is cut out by the equation $r_1r_2r_3 - r^3_4$. Thus the moduli space
of totally marked degree two rational maps is isomorphic to a Del Pezzo surface and the 
isomorphism is defined over $\mb{Z}[1/2]$.

\begin{center}
\sc{Acknowledgements}
\end{center}

The author thanks Ray Hoobler and Lloyd West for several stimulating discussions toward this paper. 
Thanks to Alon Levy for suggesting a shorter proof of Theorem 5.4. The author's research is partially 
supported by CUNY C3IRG Grant.


\end{document}